\pdfoutput=1
\RequirePackage{ifpdf}
\ifpdf 
\documentclass[pdftex]{sigma}
\else
\documentclass{sigma}
\fi

\numberwithin{equation}{section}
\newtheorem{Theorem}{Theorem}[section]
\newtheorem{Corollary}[Theorem]{Corollary}
\newtheorem{Lemma}[Theorem]{Lemma}
\newtheorem{Proposition}[Theorem]{Proposition}
{\theoremstyle{definition}
\newtheorem{Definition}[Theorem]{Definition}
\newtheorem{Example}[Theorem]{Example}
\newtheorem{Remark}[Theorem]{Remark}
\newtheorem{Notation}[Theorem]{Notation}
}

\begin{document}

\newcommand{\arXivNumber}{1403.6817}

\allowdisplaybreaks

\renewcommand{\thefootnote}{$\star$}

\renewcommand{\PaperNumber}{098}

\FirstPageHeading

\ShortArticleName{Center of Twisted Graded Hecke Algebras for Homocyclic Groups}

\ArticleName{Center of Twisted Graded Hecke Algebras\\
for Homocyclic Groups\footnote{This paper is a~contribution to the Special Issue on New Directions in Lie Theory.
The full collection is available at
\href{http://www.emis.de/journals/SIGMA/LieTheory2014.html}{http://www.emis.de/journals/SIGMA/LieTheory2014.html}}}

\Author{Wee Liang GAN~$^\dag$ and Matthew HIGHFIELD~$^\ddag$}

\AuthorNameForHeading{W.L.~Gan and M.~Highf\/ield}

\Address{$^\dag$~University of California, Riverside, CA 92521, USA}
\EmailD{\href{mailto:wlgan@math.ucr.edu}{wlgan@math.ucr.edu}}

\Address{$^\ddag$~Pepperdine University, Malibu, CA 90263, USA}
\EmailD{\href{mailto:matthew.highfield@pepperdine.edu}{matthew.highfield@pepperdine.edu}}

\ArticleDates{Received March 31, 2014, in f\/inal form October 10, 2014; Published online October 15, 2014}

\Abstract{We determine explicitly the center of the twisted graded Hecke algebras asso\-ciated to homocyclic groups.
Our results are a~generalization of formulas by M.~Douglas and B.~Fiol in [\textit{J.~High Energy Phys.} \textbf{2005} (2005), no.~9, 053, 22~pages].}


\Classification{20C08}

\renewcommand{\thefootnote}{\arabic{footnote}}
\setcounter{footnote}{0}

\section{Main results}

The notion of twisted graded Hecke algebras was introduced by S.~Witherspoon in~\cite{W}; they are variants of the
graded Hecke algebras of V.~Drinfel'd~\cite{Dr} and G.~Lusztig~\cite{Lu2} (see also~\cite{RS}) and twisted symplectic
ref\/lection algebras of T.~Chmutova~\cite{Ch}.
To a~f\/inite dimensional complex vector space~$V$, a~f\/inite subgroup~$G$ of ${\rm GL}(V)$, and a~$2$-cocycle~$\alpha$ of~$G$,
the associated twisted graded Hecke algebra $\mathsf{H}$ is, by def\/inition, a~Poincar\'e--Birkhof\/f--Witt deformation of
the crossed-product algebra $SV\#_\alpha G$, where $SV$ denotes the symmetric algebra of~$V$.
The center of $SV\#_\alpha G$ is $(SV)^G$, and it is a~natural question to determine the center of $\mathsf{H}$.
In the non-twisted case, the center of the graded Hecke algebra associated to a~f\/inite real ref\/lection group was
determined by G.~Lusztig in~\cite[Theorem~6.5]{Lu1}.
In this paper, we determine the center of $\mathsf{H}$ for the twisted graded Hecke algebra in~\cite[Example~2.16]{W},
where $V={\mathbb{C}}^n$ and~$G$ is isomorphic to a~homocyclic group $({\mathbb{Z}}/\ell{\mathbb{Z}})^{n-1}$.
(By a~homocyclic group, we mean a~direct product of cyclic groups of the same order.) In this example, the algebra
$\mathsf{H}$ is f\/initely generated as a~module over its center; the center of $\mathsf{H}$ therefore plays an important
role in the representation theory of $\mathsf{H}$.
We show that the center of $\mathsf{H}$ is generated by $n+1$ elements subject to one relation, which we determine
explicitly.
Our results are a~generalization of formulas by M.~Douglas and B.~Fiol who considered the special case when $n=3$ in
their paper~\cite{DF} on ${\mathbb{C}}^3/({\mathbb{Z}}/\ell{\mathbb{Z}})^2$ orbifolds with discrete torsion.

We state our main results in this section and give the proofs in Section~\ref{Section2}.
We shall work over ${\mathbb{C}}$.
Let~$n$ be an integer $\ge 3$, and~$\ell$ an integer $\ge 2$.
Let $V={\mathbb{C}}^n$ and let $x_1,\ldots, x_n$ be the standard basis of~$V$.
Let~$G$ be the subgroup of $SL_n({\mathbb{C}})$ consisting of all diagonal matrices~$g$ satisfying $g^\ell =1$.
Let~$\zeta$ be a~primitive~$\ell$-th root of unity.

\begin{Notation}
All subscripts are taken modulo~$n$.
For example, $x_{n+1}=x_1$.
\end{Notation}

For $i=1,\dots, n$, let $g_i$ be the element of~$G$ such that
\begin{gather*}
g_i (x_j) =
\begin{cases}
\zeta x_j, & \text{if} \quad j=i,
\\
\zeta^{-1} x_j, & \text{if} \quad j=i+1,
\\
x_j, & \mbox {else}.
\end{cases}
\end{gather*}
Observe that $g_n=g_1^{-1}\cdots g_{n-1}^{-1}$.
We have an isomorphism $({\mathbb{Z}}/\ell{\mathbb{Z}})^{n-1} \overset{\sim}{\longrightarrow} G$ def\/ined by sending
$(1,0,\dots,0),\ldots, (0,\ldots, 0,1)$ to $g_1,\dots, g_{n-1}$, respectively.

Def\/ine the 2-cocycle $\alpha: G\times G\to {\mathbb{C}}^\times$ of~$G$~by
\begin{gather*}
\alpha\big(g_1^{i_1} \cdots g_{n-1}^{i_{n-1}}, g_1^{j_1} \cdots g_{n-1}^{j_{n-1}}\big) = \zeta^{-i_1j_2-i_2j_3-\dots -
i_{n-2}j_{n-1}}.
\end{gather*}
If~$E$ is an algebra, an action of~$G$ on~$E$ is a~homomorphism $G\to \mathrm{Aut}(E)$.
Recall that for any algebra~$E$ and an action of~$G$ on~$E$, one has the crossed product algebra $E\#_\alpha G$.
As a~vector space, $E\#_\alpha G$ is $E\otimes{\mathbb{C}} G$; the product is def\/ined~by
\begin{gather*}
(r\otimes g)(s\otimes h) = \alpha (g,h) r (g\cdot s) \otimes gh
\end{gather*}
for all $r,s \in E$ and $g,h\in G$.
If $g, h\in G$, then we shall denote their product in $E\#_\alpha G$ by $g*h$; thus,
\begin{gather*}
g* h = \alpha(g,h) gh.
\end{gather*}
One has, for any $i, j \in \{1,\ldots,n\}$ with $|i-j|\notin\{1, n-1\}$,
\begin{gather*}
g_{i+1}*g_i = \zeta g_i * g_{i+1},
\qquad
g_i *g_j = g_j *g_i.
\end{gather*}
Let $t=(t_1,\dots, t_n)\in {\mathbb{C}}^n$, and write $TV$ for the tensor algebra of~$V$.
Following~\cite[Examp\-le~2.16]{W}, we make the following def\/inition.

\begin{Definition}
\label{hhh}
Let $\mathsf{H}$ be the associative algebra def\/ined as the quotient of $TV\#_\alpha G$ by the relations:
\begin{gather*}
x_ix_{i+1}- x_{i+1} x_i = t_i g_i,
\qquad
x_ix_j - x_jx_i = 0
\end{gather*}
for all $i, j \in \{1,\ldots,n\}$ with $|i-j|\notin\{1, n-1\}$.
\end{Definition}

\begin{Remark}
By~\cite[Theorem 2.10]{W} and~\cite[Example 2.16]{W}, the algebra $\mathsf{H}$ in Def\/inition~\ref{hhh} is a~twisted
graded Hecke algebra for~$G$.
(However, when $n>3$ and $\ell=2$, this is not the most general twisted graded Hecke algebra for~$G$; see~\cite[Example 2.16]{W}
and~\cite[Example 5.1]{W2}.)
\end{Remark}

Let ${\mathbb{C}}[y_1^\pm, \ldots, y_n^\pm]$ be the algebra of Laurent polynomials in the variables $y_1,\ldots, y_n$.
The group~$G$ acts on ${\mathbb{C}}[y_1^\pm, \ldots, y_n^\pm]$~by
\begin{gather*}
g_i y_1^{p_1}\cdots y_n^{p_n} = \zeta^{p_i-p_{i+1}}y_1^{p_1}\cdots y_n^{p_n}
\end{gather*}
for all $i\in\{1,\ldots, n-1\}$ and $p_1,\ldots,p_n\in {\mathbb{Z}}$.

\begin{Proposition}
\label{T1}
There is an injective homomorphism
\begin{gather*}
\Theta: \ \mathsf{H} \longrightarrow {\mathbb{C}}[y_1^\pm, \ldots, y_n^\pm]\#_\alpha G
\end{gather*}
such that
\begin{gather}
\Theta(x_i) = y_i - \left(\frac{\zeta t_i}{\zeta-1} \right) y_{i+1}^{-1}g_i,
\label{theta1}
\\
\Theta(g_i) = g_i
\label{theta3}
\end{gather}
for all $i\in\{1,\ldots, n\}$.
\end{Proposition}

Let
\begin{gather*}
I = \{\{i_1<\dots<i_k\} \mid k\ge 0; \ i_1,\ldots, i_k\in \{1,\dots,n\} \},
\\
J = \{\{i_1<\dots<i_k\}\in I \mid |i_r-i_s|\notin\{1,n-1\}\ \text{for all}\ r, s \}.
\end{gather*}
Def\/ine the elements $\delta$, $\varepsilon_1, \ldots, \varepsilon_n$ of ${\mathbb{Z}}^n$~by
\begin{gather*}
\delta = (1,1,\dots,1),
\quad
\varepsilon_1 = (1,1,0,\dots,0),
\quad
\varepsilon_2 = (0,1,1,0,\dots),
\quad
\ldots,
\quad
\varepsilon_n = (1,0,\dots,0,1).
\end{gather*}

\begin{Notation}
For any variables $\omega_1, \ldots, \omega_n$ and $p=(p_1,\ldots,p_n)\in {\mathbb{Z}}^n$, we denote by $\omega^p$ the
expression $\omega_1^{p_1} \cdots \omega_n^{p_n}$.
\end{Notation}

We shall set
\begin{gather*}
\tau_i = \frac{t_i}{\zeta-1}
\quad
\text{for}
\quad
i=1,\ldots, n-1,
\qquad
\tau_n = \frac{\zeta t_n}{\zeta-1}.
\end{gather*}
Def\/ine the element $w \in \mathsf{H}$~by
\begin{gather*}
w = \sum\limits_{\{i_1<\dots<i_k\}\in J} \tau_{i_1}\cdots \tau_{i_k} x^{\delta-\varepsilon_{i_1} - \dots -
\varepsilon_{i_k}} g_{i_1}*\cdots * g_{i_k}.
\end{gather*}

\begin{Example}
If $n=3$, then
\begin{gather*}
w = x_1x_2x_3 +\tau_1 x_3 g_1 + \tau_2 x_1 g_2 + \tau_3 x_2 g_3
= x_1 x_2 x_3 +\frac{1}{\zeta-1} \left(t_1 x_3 g_1 + t_2 x_1 g_2 + \zeta t_3 x_2 g_3 \right).
\end{gather*}
In particular, if $n=3$ and $\ell=2$, the formula for~$w$ is in~\cite[Lemma~7.1]{CGW}.
\end{Example}

\begin{Theorem}
\label{T2}
The center of $\mathsf{H}$ is generated as an algebra by $x_1^\ell, \ldots, x_n^\ell$, and~$w$.
\end{Theorem}

Let $\mathsf{Z}$ be the center of $\mathsf{H}$.
For $ r=0,\ldots, \lfloor \ell/2\rfloor$, set
\begin{gather*}
\nu_r = (-1)^r \frac{\ell}{\ell-r} \binom{\ell-r}{r},
\end{gather*}
and set
\begin{gather*}
{\widetilde{\tau}}_i = \tau_i^\ell
\quad \text{for} \quad
i=1, \ldots, n-1,
\qquad
{\widetilde{\tau}}_n = (-1)^{n(\ell-1)}\tau_n^\ell.
\end{gather*}

We def\/ine a~polynomial~$F$ in the $n+1$ variables $a_1, \ldots, a_n$ and~$b$~by
\begin{gather}
\label{relation}
F = \sum\limits_{\{i_1<\dots<i_k\}\in J} {\widetilde{\tau}}_{i_1}\cdots {\widetilde{\tau}}_{i_k}
a^{\delta-\varepsilon_{i_1} - \dots - \varepsilon_{i_k}} - \sum\limits_{r=0}^{\lfloor \ell/2 \rfloor} (-1)^{nr}
\zeta^{(n-2)r} \nu_r (\tau_1\cdots \tau_n)^r b^{\ell-2r}.
\end{gather}

\begin{Corollary}
\label{maincor}
The assignment
\begin{gather}
\label{mapphi}
a_i\mapsto x_i^\ell
\quad \text{for} \quad
i=1, \ldots, n,
\qquad
b \mapsto w
\end{gather}
defines an isomorphism
\begin{gather}
\label{center}
{\mathbb{C}}[a_1,\ldots, a_n,b]/(F) \stackrel{\sim}{\longrightarrow} \mathsf{Z}.
\end{gather}
\end{Corollary}

In the undeformed case, when $t_1=\dots=t_n=0$, the polynomial~$F$ is equal to $a_1\cdots a_n - b^\ell$.

\section{Proof of main results}\label{Section2}

\begin{proof}
[Proof of Proposition~\ref{T1}]
For $i=1,\ldots, n-1$, we def\/ine $\Theta(x_i)$, $\Theta(x_n)$, and $\Theta(g_i)$
by~\eqref{theta1} and~\eqref{theta3}.
It follows from a~straightforward verif\/ication that~$\Theta$ is a~well-def\/ined homomorphism.

It remains to see that~$\Theta$ is injective.
Observe that $\mathsf{H}$ is spanned by the monomials $x^p g$ for $p=(p_1,\dots,p_n)\in {\mathbb{Z}}^n$ and $g\in G$,
where $p_1,\ldots, p_n\ge 0$.
We call $p_1+\dots+p_n$ the total degree of the monomial $x^p g$.
The image of $x^p g$ under~$\Theta$ is the sum of $y^p g$ with terms of strictly smaller total degrees.
Therefore, if $\alpha\in \mathsf{H}$ is nonzero, we can write it as a~sum $\alpha_0 + \alpha_1 + \cdots$, where
$\alpha_k$ is a~linear combination of monomials~$x^p g$ with total degree~$k$.
If~$k$ is the maximal integer with~$\alpha_k$ nonzero, then $\Theta(\alpha_k)$ is nonzero, and hence $\Theta(\alpha)$ is
also nonzero.
\end{proof}

\begin{Remark}
\label{pbwbasis}
It follows from Proposition~\ref{T1} that the monomials $x_1^{p_1}\cdots x_n^{p_n} g$ for non-negative integers $p_1,
\ldots, p_n$ and $g\in G$ form a~basis for $\mathsf{H}$ (called the PBW basis of $\mathsf{H}$).
This was f\/irst proved in~\cite[Example 2.16]{W} using~\cite[Theorem 2.10]{W}.
\end{Remark}

We have an increasing f\/iltration on $\mathsf{H}$ def\/ined by setting $\deg(x_i)=1$ and $\deg(g)=0$ for all
$i\in\{1,\ldots,n\}$, $g\in G$.
It is immediate from Remark~\ref{pbwbasis} that the natural homomorphism $SV\#_\alpha G \to {\mathrm{gr}} \mathsf{H}$ is
an isomorphism, where ${\mathrm{gr}} \mathsf{H}$ denotes the associated graded algebra of $\mathsf{H}$.

The proof of~\eqref{thetaw} in the following lemma is the key calculation in this paper.

\begin{Lemma}\label{keylemma}\quad
\begin{enumerate}\itemsep=0pt
\item[$(i)$] One has:
\begin{gather}
\Theta\big(x_i^\ell\big) = y_i^\ell - \tau_i^\ell y_{i+1}^{-\ell},
\label{thetaxil}
\\
\Theta\big(x_n^\ell\big) = y_n^\ell - (-1)^{n(\ell-1)} \tau_n^\ell y_1^{-\ell},
\label{thetaxnl}
\end{gather}
for all $i\in \{1,\ldots, n-1\}$.

\item[$(ii)$] One has:
\begin{gather}
\label{thetaw}
\Theta(w) = y_1\cdots y_n + (-1)^n \zeta^{n-2} \tau_1\cdots \tau_n y_1^{-1} \cdots y_n^{-1}.
\end{gather}
\end{enumerate}
\end{Lemma}

\begin{proof}
(i) To prove~\eqref{thetaxil}, we need to show that
\begin{gather}
\label{eqexpand1}
\underbrace{\big(y_i - \zeta \tau_i y_{i+1}^{-1}g_i\big)
\cdots
\big(y_i - \zeta \tau_i y_{i+1}^{-1}g_i\big)}_\ell = y_i^\ell - \tau_i^\ell y_{i+1}^{-\ell}.
\end{gather}
Since $g_iy_i=\zeta y_i g_i$ and $g_i y_{i+1}^{-1} = \zeta y_{i+1}^{-1} g_i$, the product on the left hand side
of~\eqref{eqexpand1} is a~linear combination of $y_i^k y_{i+1}^{k-\ell} g_i^{\ell-k}$ for $k=0, 1,\ldots, \ell$.
Moreover, the coef\/f\/icient of $y_i^k y_{i+1}^{k-\ell} g_i^{\ell-k}$ in this linear combination is the same as the
coef\/f\/icient of $u^k$ when we expand the product
\begin{gather}
\label{poly1}
\big(u- \zeta^\ell \tau_i\big) \big(u- \zeta^{\ell-1} \tau_i\big) \cdots (u- \zeta \tau_i)
\end{gather}
in the polynomial ring ${\mathbb{C}}[u]$.
Since the polynomial in~\eqref{poly1} is equal to $u^\ell-\tau_i^\ell$, the identity~\eqref{thetaxil} follows.
The proof of~\eqref{thetaxnl} is similar except that
\begin{gather*}
\underbrace{g_n * \dots * g_n}_\ell = (-1)^{n(\ell-1)}.
\end{gather*}
(ii) For any $h_*=\{h_1<\dots<h_j\} \in I$, we let
\begin{gather*}
h_*' = \{h_r\in h_* \mid h_s-h_r \in \{1, 1-n\}\ \text{for some}\ s \},
\\
\chi(h_*) = |\{h_r\in h'_* \mid h_r\neq n \}| - |\{h_r\in h'_* \mid h_r=n \}|,
\\
E(h_*) = \zeta^{\chi (h_*)} \tau_{h_1} \cdots \tau_{h_j} y^{\delta - \varepsilon_{h_1} - \dots - \varepsilon_{h_j}}
g_{h_1} * \cdots * g_{h_j}.
\end{gather*}
Now suppose $i_*=\{i_1<\dots<i_k\}\in J$.
Let~$D$ be the subset of $\{1,\dots,n\}$ consisting of all~$d$ such that $d \not\equiv i_r,\ i_r+1$ (mod~$n$) for all~$r$.
We denote by $d_1 < \dots < d_p$ the elements of~$D$.
Then
\begin{gather*}
\Theta\big(\tau_{i_1}\cdots \tau_{i_k} x^{\delta-\varepsilon_{i_1} - \dots - \varepsilon_{i_k}} g_{i_1}*\cdots * g_{i_k}\big)
\\
\qquad
= \tau_{i_1}\cdots \tau_{i_k} \left(y_{d_1} - \frac{\zeta t_{d_1}}{\zeta-1}y_{d_1+1}^{-1}g_{d_1} \right) \cdots
\left(y_{d_p} - \frac{\zeta t_{d_p}}{\zeta-1}y_{d_p+1}^{-1}g_{d_p} \right) g_{i_1}*\cdots * g_{i_k}
\\
\qquad
= \tau_{i_1}\cdots \tau_{i_k} \sum\limits_{S\subset D} Y_{d_1}(S) \cdots Y_{d_p}(S) g_{i_1}*\cdots * g_{i_k},
\end{gather*}
where, for $r=1,\ldots, p$,
\begin{gather*}
Y_{d_r}(S) =
\begin{cases}
y_{d_r}, & \text{if} \quad d_r\notin S,
\\
- \zeta(\zeta-1)^{-1} t_{d_r}y_{d_r+1}^{-1}g_{d_r}, & \text{if} \quad d_r\in S.
\end{cases}
\end{gather*}
Setting $h_* = i_* \cup S$, we obtain\footnote{Note that if $d_r\in S$ but $d_r+1\in D-S$, then the term $g_{d_r}$ in
$Y_{d_r}(S)$ appears on the left of the term $y_{d_r+1}$ of $Y_{d_r+1}(S)$ and one has $g_{d_r} y_{d_r+1} = \zeta^{-1}
y_{d_r+1} g_{d_r}$.
However, if $n\in S$ but $1\in D-S$, then the term $g_n$ in $Y_n(S)$ already appears to the right of the term $y_1$ of
$Y_1(S)$.
This is the reason why the def\/inition of $\tau_n$ dif\/fers from the corresponding def\/initions of $\tau_1,\ldots,
\tau_{n-1}$ by a~factor of~$\zeta$.}
\begin{gather*}
\Theta\big(\tau_{i_1}\cdots \tau_{i_k} x^{\delta-\varepsilon_{i_1} - \dots - \varepsilon_{i_k}} g_{i_1}*\cdots * g_{i_k}\big) =
\sum\limits_{\{h_*\in I \mid i_*\subset h_*- h'_* \}} (-1)^{|h_*|-|i_*|} E(h_*).
\end{gather*}
Hence,
\begin{gather*}
\Theta(w) = \sum\limits_{\{i_1<\dots<i_k\}\in J} \Theta\big(\tau_{i_1}\cdots \tau_{i_k} x^{\delta-\varepsilon_{i_1} - \dots -
\varepsilon_{i_k}} g_{i_1}*\dots * g_{i_k}\big)
\\
\hphantom{\Theta(w)}{}
{} = \sum\limits_{i_*\in J} \left(\sum\limits_{\{h_*\in I \mid i_*\subset h_*- h'_*\}} (-1)^{|h_*|-|i_*|} E(h_*) \right)
 = \sum\limits_{h_*\in I} \left (E(h_*) \sum\limits_{i_*\subset h_*-h'_*} (-1)^{|h_*|-|i_*|} \right).
\end{gather*}
If $|h_*|=n$, then $h'_*= h_*$.
If $|h_*|\notin\{0,n\}$, then $h'_* \neq h_*$.
Therefore,
\begin{gather*}
E(h_*) \sum\limits_{i_*\subset h_*-h'_*} (-1)^{|h_*|-|i_*|} =
\begin{cases}
y_1\cdots y_n & \text{if} \quad |h_*|=0,
\\
(-1)^n \zeta^{n-2} \tau_1\cdots \tau_n y_1^{-1} \cdots y_n^{-1} & \text{if} \quad |h_*|=n,
\\
0 & \text{else.}
\end{cases}
\tag*{\qed}
\end{gather*}
\renewcommand{\qed}{}
\end{proof}

\begin{proof}[Proof of Theorem~\ref{T2}] It is easy to see that the center of $SV\#_\alpha G$ is the algebra of~$G$-invariant
elements $(SV)^G$ of $SV$, and moreover, the algebra $(SV)^G$ is generated by $x_i^\ell$ ($i=1,\ldots,n$) and $x_1\cdots x_n$.

Using Lemma~\ref{keylemma}, we see that
\begin{gather*}
\Theta\big(x_i^\ell\big)
\quad
\text{for}\quad i=1,\ldots, n,
\qquad
\text{and}
\qquad
\Theta(w)
\end{gather*}
are in the center of ${\mathbb{C}}[y_1^\pm, \ldots, y_n^\pm]\#_\alpha G$.
Since the homomorphism~$\Theta$ is injective, the ele\-ments~$x_i^\ell$ ($i=1,\ldots, n$) and~$w$ are in the center of
$\mathsf{H}$.
Since the principal symbols of $x_1^\ell, \ldots, x_n^\ell$ and~$w$ in $SV\#_\alpha G$ are, respectively, $x_1^\ell,
\ldots, x_n^\ell$ and $x_1\cdots x_n$, the theorem follows from a~standard argument.
\end{proof}

\begin{proof}
[Proof of Corollary~\ref{maincor}] Let ${\widetilde{a}}_1 = \Theta(x_1^\ell)$, $\ldots$, ${\widetilde{a}}_n =
\Theta(x_n^\ell)$, and ${\widetilde{b}}=\Theta(w)$.
By Lemma~\ref{keylemma},
\begin{gather*}
{\widetilde{a}}_i = y_i^\ell - {\widetilde{\tau}}_i y_{i+1}^{-\ell}
\quad
\text{for}
\quad
i=1,\ldots, n,
\\
{\widetilde{b}} = y_1\cdots y_n + (-1)^n \zeta^{n-2} \tau_1\cdots \tau_n y_1^{-1} \cdots y_n^{-1}.
\end{gather*}
By a~calculation completely similar to the proof of~\eqref{thetaw}, one has
\begin{gather}
\label{leftside}
\sum\limits_{\{i_1<\dots<i_k\}\in J} {\widetilde{\tau}}_{i_1}\cdots {\widetilde{\tau}}_{i_k}
{\widetilde{a}}^{\delta-\varepsilon_{i_1} - \dots - \varepsilon_{i_k}} = (y_1 \cdots y_n)^\ell + (-1)^{n\ell}
(\tau_1\cdots \tau_n)^\ell (y_1\cdots y_n)^{-\ell}.
\end{gather}
We claim that we also have
\begin{gather}
\label{rightside}
\sum\limits_{r=0}^{\lfloor \ell/2 \rfloor} (-1)^{nr} \zeta^{(n-2)r} \nu_r (\tau_1\cdots \tau_n)^r
{\widetilde{b}}^{\ell-2r} = (y_1 \cdots y_n)^\ell + (-1)^{n\ell} (\tau_1\cdots \tau_n)^\ell (y_1\cdots y_n)^{-\ell}.
\end{gather}
To see this, recall that the Chebyshev polynomials of the f\/irst kind are def\/ined recursively by $T_0(\xi) = 1$,
$T_1(\xi) = \xi$, and
\begin{gather*}
T_m(\xi) = 2\xi T_{m-1}(\xi) - T_{m-2}(\xi)
\quad
\text{for}
\quad
m=2, 3, \ldots.
\end{gather*}
It is well known (and can be easily proved by induction) that
\begin{gather}
2 T_\ell \left(\frac{\xi}{2}\right) = \sum\limits_{r=0}^{\lfloor \ell/2 \rfloor} \nu_r \xi^{\ell-2r},
\label{che1}
\\
2T_\ell \left(\frac{\xi+ \xi^{-1}}{2}\right) = \xi^\ell + \xi^{-\ell}.
\label{che2}
\end{gather}
By~\eqref{che1} and~\eqref{che2}, one has the identity
\begin{gather*}
\xi^\ell + \xi^{-\ell} = \sum\limits_{r=0}^{\lfloor \ell/2 \rfloor} \nu_r \big(\xi+\xi^{-1}\big)^{\ell-2r},
\end{gather*}
and hence the identity
\begin{gather*}
\xi^\ell + \varrho^{2\ell} \xi^{-\ell} = \sum\limits_{r=0}^{\lfloor \ell/2 \rfloor} \nu_r \varrho^{2r} \big(\xi+\varrho^2
\xi^{-1}\big)^{\ell-2r}
\end{gather*}
where $\xi$ and~$\varrho$ are formal variables.
By setting $\xi=y_1\cdots y_n$ and choosing~$\varrho$ to be a~square-root of $(-1)^n \zeta^{n-2} \tau_1\cdots \tau_n$, we
obtain~\eqref{rightside}.

By Proposition~\ref{T1}, Theorem~\ref{T2}, and the equations~\eqref{leftside} and~\eqref{rightside}, the
assignment~\eqref{mapphi} def\/ines a~surjective homomorphism
\begin{gather*}
\Phi: \  {\mathbb{C}}[a_1,\ldots, a_n,b] \to \mathsf{Z}
\end{gather*}
such that $\Phi(F)=0$.
Suppose $D\in {\mathbb{C}}[a_1,\ldots, a_n,b]$ and $\Phi(D)=0$.
We can write
\begin{gather*}
D = \sum\limits_{r=0}^{\ell-1} D_r(a_1,\dots,a_n) b^r + R,
\end{gather*}
where $D_r (a_1,\dots,a_n) \in {\mathbb{C}}[a_1,\ldots, a_n]$ for $r=0,\ldots, \ell-1$, and $R\in (F)$.
Thus,
\begin{gather}
\label{eqnd}
\sum\limits_{r=0}^{\ell-1} D_r \big(x_1^\ell,\dots, x_n^\ell\big) w^r = 0.
\end{gather}
We claim that $D_r(a_1,\dots,a_n)=0$ for all~$r$.
Suppose not; then let~$m$ be the maximal integer such that $D_m(a_1,\dots,a_n)\neq 0$.
Let $x_1^{\ell p_1} \cdots x_n^{\ell p_n}$ be a~monomial in $D_m(x_1^\ell,\dots, x_n^\ell)$ with nonzero coef\/f\/icient.
Since $0\le m<\ell$, when we write the left hand side of~\eqref{eqnd} in terms of the PBW basis, the coef\/f\/icient of
$x_1^{\ell p_1+m} \cdots x_n^{\ell p_n+m}$ is nonzero, a~contradiction.
Hence, the kernel of~$\Phi$ is~$(F)$.
This proves $\eqref{center}$.
\end{proof}

\begin{Remark}
When $n=3$, the algebra $\mathsf{H}$ is Morita equivalent to a~deformed Sklyanin algebra~$S_{\rm def}$ def\/ined by C.~Walton
in~\cite[Def\/inition~IV.2]{Wa}.
More precisely, if $n=3$ and
\begin{gather*}
e= \frac{1}{\ell} \sum\limits_{r=0}^{\ell-1}g_1^r,
\end{gather*}
one has $\mathsf{H}e\mathsf{H} = \mathsf{H}$ and $e\mathsf{H}e \cong S_{\rm def}$ where the parameters for $S_{\rm def}$
(following the notations in~\cite[Def\/inition~IV.2]{Wa}) are $a=1$, $b=\zeta$, $c=d_i=0$, and $e_i=-\zeta t_i$ for $i=1,2,3$.
This follows from the observation that, for $n=3$, setting $\phi_i = x_i g_{i+1}$, one has $\phi_i\phi_{i+1} - \zeta
\phi_{i+1}\phi_i = \zeta t_i$ for all~$i$.
The algebra $S_{\rm def}$ (with above parameters) was f\/irst studied by M.~Douglas and B.~Fiol, see~\cite[(3.10)]{DF}.
Our formulas~\eqref{theta1}--\eqref{theta3} are a~generalization of~\cite[(4.6)]{DF}, and our equation~\eqref{relation}
is a~generalization of~\cite[(4.7)]{DF}.
The formulas in~\eqref{thetaxil}--\eqref{thetaw} are generalizations of~\cite[(4.8)]{DF}.
\end{Remark}

\subsection*{Acknowledgements}

We thank the referees for their many helpful comments.

\pdfbookmark[1]{References}{ref}
\LastPageEnding

\end{document}